\documentclass [11pt]{article}
\newfont{\bbb} {msbm10}
\newcommand{\Bbb}[1]{\mbox{\bbb#1}}
\newcommand{\R}{\Bbb{R}}
\newcommand{\bS}{\Bbb{S}}

\newcommand{\Z}{\Bbb{Z}}

\newcommand{\D}{\Bbb{D}}
\newcommand{\C}{\Bbb{C}}

\newcommand{\sm}{\setminus}
\newcommand{\sbs}{\subset}
\newcommand{\ra}{\rightarrow}

\newcommand{\tQ}{\tilde{Q}}
\newcommand{\ctQ}{\overline {(\tilde{Q})}}

\newcommand{\tS}{\tilde{S}}

\newcommand{\tN}{\tilde{N}}
\newcommand{\tA}{\tilde{A}}
\newcommand{\tW}{\tilde{W}}

\newcommand{\p}{\partial}
\newcommand{\bo}{\partial_\infty}

\newcommand{\met}{{\cal{MET}}}

\newcommand{\mo}{{\cal{MET}}^{\, sec\, < \, 0}}
\newcommand{\moo}{{\cal{MET}}^{\, sec\, \leq \, 0}}

\newcommand{\cL}{{\cal{L}}}
\newcommand{\cC}{{\cal{C}}}

\newcommand{\cA}{{\cal{A}}}
\newcommand{\cB}{{\cal{B}}}

\usepackage[all]{xy}

\begin{document}

\title{  The Space of Nonpositively Curved Metrics of a Negatively Curved Manifold}
\author{F. T. Farrell and P. Ontaneda\thanks{Both authors were
partially supported by NSF grants.}}
\date{}

\maketitle

\begin{abstract} We show that the space of nonpositively curved metrics of a closed
negatively curved Riemannian $n$-manifold, $n\geq 10$, is highly
non-connected. 
\end{abstract}
\vspace{.3in}

\noindent {\bf \Large  Section 0. Introduction.}\\

Let $M^n$ be a  closed smooth manifold of dimension $dim\, M =n$.
We denote by $\met (M)$ the space of all smooth Riemannian metrics
on $M$ and we consider $\met (M)$ with the smooth topology. Also, we
denote by $\mo(M)$ the subspace formed by
all negatively curved Riemannian
metrics on $M$. In \cite{FO6} we proved that $\mo (M)$
always has infinitely many path-components,  provided $n\geq 10$
and it is non-empty.
Moreover we showed that all the groups $\pi_{2p-4}(\mo (M))$ are non-trivial
for every prime number $p>2$, and such that $p<\frac{n+5}{6}$
(this is true in every component of $\mo(M)$).
In fact, these groups contain the infinite sum $(\Z_p)^\infty$ of $\Z_p=\Z/p\,\Z$'s.
We also showed that $\pi_1(\mo (M))$ contains the infinite sum $(\Z_2)^\infty$ when $n\geq 12$ (see also \cite{FO7}).
All these results follow from the Main Theorem in \cite{FO6}, which states
that the orbit map $\Lambda_g:DIFF(M)\ra \mo(M)$ is ``very non-trivial"
at the $\pi_k$-level. Here $DIFF(M)$ is the group of self-diffeomorphims
on $M$ and $\Lambda_g(\phi)=\phi_*g$ (see the introduction of \cite{FO6}
for more details).\\

Let $\moo(M) $ be the subspace of $\met(M)$ formed by all non-positively curved Riemannian metrics on $M$.
In this paper we generalize to $\moo(M)$ the results mentioned above, provided $\pi_1M$ is (word) hyperbolic:\\

\noindent {\bf Main Theorem.} {\it Let $M^n$ be a closed smooth manifold with hyperbolic fundamental group $\pi_1 M$.
Assume $\moo(M)$ is non-empty. Then 
\begin{enumerate}
\item[(i)] the space $\moo(M)$ has infinitely many components, 
provided $n\geq 10$.
\item[(ii)] The  group $\pi_1(\moo (M^n))$ is not trivial
when $n\geq 12$. In fact it contains the infinite sum $(\Z_2)^\infty$ as a subgroup. 
\item[(iii)] The groups $\pi_{2p-4}(\moo (M^n))$ are non-trivial
for every prime number $p>2$, and such that $p<\frac{n+5}{6}$.
In fact, these groups contain the infinite sum $(\Z_p)^\infty$ as a subgroup.
\end{enumerate}}
\newpage

\noindent{\bf Remarks.}

%\noindent {\bf 1.} The result in \cite{FO6} states, for the case $k=1$, that  %$n\geq 14$ is needed. But this can be improved to
%$n\geq 12$, see \cite{FO7}.

\noindent {\bf 1.} The results for $\pi_k\,\moo\,(M)$, $k>0$, given above are true
relative to any base point, that is, for every component of $\moo(M)$.
%(The space $\moo\,(M)$ is locally path-connected, hence
%connected  components coincide with path-connected components. TRUE?)

\noindent {\bf 2.} The decoration ``$sec\, \leq \, 0$'' can be tightened to ``$a\leq sec\leq\, 0$'', for any $a< 0$.

\noindent {\bf 3.} The Theorem above follows from a nonpositively curved version of the
Main Theorem of \cite{FO6} (which we do not state to save space).
This nonpositively curved version is obtained from the Main Theorem in \cite{FO6} by
replacing $\mo(M)$ by $\moo(M)$ and adding the hypothesis
``$\pi_1(M)$ is hyperbolic". This is the result that we prove in this
paper. And, as in \cite{FO6}, we obtain the following corollary.\\

\noindent {\bf Corollary.} {\it
Let $M$ be a closed smooth $n$-manifold with  
$\pi_1(M)$ hyperbolic. 
Let $I\sbs(-\infty, 0]$ and assume that $\met^{sec\in I}(M)$ is not empty.
Then the inclusion map  $\met^{sec\in I}(M)\hookrightarrow\moo(M)$ 
is not null-homotopic, provided $n\geq 10$. 

Moreover, the induced maps of
this inclusion, at the $k$-homotopy level, are not constant for $k=0$, and non-zero for $k$ and $n$ as in
cases (ii.), (iii.) in the Main Theorem. Furthermore, the image of these
maps satisfy a statement analogous to the one in the Addendum  to the Main Theorem in \cite{FO6}.}\\

Here $\met^{sec\in I}$ has the obvious meaning.
In particular taking $I=(-\infty,0)$ we get that the inclusion
$\mo(M)\hookrightarrow\moo(M)$ is not nullhomotopic,
provided  $n\geq 10$
and $M$ admits a negatively curved metric.\\ \\

In some sense it is quite surprising that we were able to extend the results in \cite{FO6} to the nonpositively curved case because negative curvature
is a ``stable" condition (the space $\mo(M)$ is open in $\met(M)$) while 
$\moo$ is not stable. Indeed it is not even known whether $\moo(M)$ is locally
contractible or even locally connected. We state these as questions:\vspace{.2in}

\noindent {\bf Questions.}\\
{\bf 1.} Is the space $\moo(M)$ of nonpositively curved metrics on $M$
locally contractible?\\
{\bf 2.} Is the space $\moo(M)$ of nonpositively curved metrics on $M$
locally connected?\vspace{.2in}

So, we prove here that $\moo(M^n)$ is not (globally) connected
when it is not empty, $n\geq 10$ and $\pi_1M$ is hyperbolic. But on the other hand it is not known whether $\moo(M)$ is locally connected.\\

In this paper there are two additional obstacles to pass from negative curvature to nonpositive curvature.
First, since we can now have parallel geodesic rays emanating perpendicularly from a closed
geodesic, the obstructions we defined in \cite{FO6} (which lie in the pseudoisotopy space of
$\bS^\times\bS^{n-2}$) may not be homeomorphisms at infinity.\\

The second problem is that we may now have a whole family of closed
geodesics freely homotopic to a given one. But in our previous papers we strongly
used the fact there is a {\it unique} such closed geodesic. Moreover, we strongly used the fact
that such unique closed geodesics depend smoothly on the metric. This does not happen in nonpositive
curvature. Even worse: there are examples of smooth families $g_t$, $t\in [0,1]$, of nonpositively curved metrics
such that there is no continuous path of closed $g_t$-geodesics joining a closed $g_1$-geodesic
to a closed $g_0$-geodesic (all closed geodesics in the same free homotopy class). See for instance the ``swinging neck" in Appendix  A. We deal with this by incorporating the closed geodesics
into the system, but we pay a price for this: instead of dealing with discs (to prove that an element is zero in a homotopy
group) we have to deal with  more complicated spaces which we call ``cellular discs''. Because of this the use of shape theory becomes necessary.\\

In section 1 we define cellular discs and give some preliminary results.
In section 2 we prove the Main Theorem (see remark 3 above).
%In section 4 we prove the second part of Theorem B.
We shall refer to \cite{FO6} for some details.\\

We are grateful to Ross Geoghegan and Jerzy Dydak for the valuable information provided to us.\vspace{.6in}

\noindent {\bf \Large  Section 1. Preliminaries.}\\ \\

\noindent {\bf \large A.  Cellular Discs.}\\

We will consider the $k$-disc $\D^k=\{ \, x\in \R^k\,\, :\,\, |x|\leq 1 \}$ with base point $u_0=(1,0,0,...,0)$.
A {\it cellular $k$-disc} is a metrizable compact pointed topological space $(X, x_0)$ together with a surjective continuous map
$\eta:(X, x_0)\ra(\D^k, u_0)$ such that the pre-image $\eta^{-1}(u)$, $u\in \D^k$, is homeomorphic to the $\ell_u$-disc $\D^{\ell_u}$,
with $0\leq\ell_u\leq \ell$, for some $\ell<\infty$ and all $u\in\D^k$. 
We write $X_u=\eta^{-1}(u)$, $X_0=\eta^{-1}(u_0)$ and $\p X=\eta^{-1}(\p\D^k)=\eta^{-1}(\bS^{k-1})$.\\
%We write $\p X=\eta^{-1} (\bS^{k-1})$, $X_u=\eta^{-1}(u)$ and $X_0=X_{x_0}$.\\

A pair $((X,x_0),X')$, $x_0\in X'\sbs \p X$, together with a map $\eta:X\ra \D^k$ is a {\it cellular $k$-disc pair} if
$X$ (that is $((X,x_0),\eta)$) is a cellular $k$-disc and  $\p X$ is fibered homeomorphic to
$X'\times X_0$, that is, there is a homeomorphism $X'\times X_0\ra\p X$ that sends $\{x'\}\times X_0$ to $X_u$, $u=\eta(x')\in \bS^{k-1}$. In particular $\eta|_{X'}:X'\ra\bS^{k-1}$ is a homeomorphism.
We identify $X'$ with $\bS^{k-1}$ and say that $(X,\bS^{k-1})$ is a cellular $k$-disc pair.\\

%Sometimes we will say that $X$ is a cellular disc (instead of the pair
%$((X,x_0),\eta)$). If $X$ satisfies only condition (i) above we say
%that $X$ is a {\it weakly cellular disc}.
%We identify $\p X$ with $\bS^{k-1}\times X_0$ and $\bS^{k-1}\times\{ x_0\}\sbs\p X$ with $\bS^{k-1}$. 
Note that it makes sense to say that a map $h:\bS^{k-1}\ra Y$ extends to a cellular
$k$-disc pair $(X,\bS^{k-1})$.\\

In the proofs of the following two Propositions we use shape theory (see for instance
\cite{Dydak-Segal}, \cite{Mardesic-Segal}). Recall that the objects of the shape category are
pointed spaces, and for two such objects $A$ and  $B$ we denote the set of morphisms by $sh\{ A, B\}$.  There is a 
functor, the {\it shape functor}, from the pointed homotopy category of topological spaces to the shape category. Hence, for each pair of
pointed spaces $A$ and $B$ we get a {\it shape map}  between  $[A,B]$, the set of
pointed homotopy classes of maps, and  $sh\{ A,B\}$. In particular there are shape maps from the homotopy groups
of $B$ to the homotopy pro-groups of $B$ (these are the shape versions of the homotopy groups of $B$).\\

Recall that a metric space Z is 
%{\it uniformly weakly locally contractible}
$LC^m$
 if for every $z\in Z$ and $\epsilon>0$ there is a $\delta>0$ 
such that any continuous map $f:P\ra B_\delta(z)$, $P$ a locally finite polyhedron of dimension $\leq m$, is homotopic in $B_\epsilon(z)$ to a constant map. And $Z$ is $LC^\infty$ if it is $LC^m$ for every $m$. 
We will use the following facts:\\

\noindent {\bf Fact 1.} {\it A cell-like map between finite dimensional spaces is a shape equivalence}
\cite{Sher}.\\

\noindent {\bf Fact 2.} {\it Let $W$ and $Z$ be pointed spaces. Assume 
$W$ is finite dimensional and $Z$ is sufficiently nice (for instance $Z$ is $LC^\infty$). 
Then} $$ [W,Z]\stackrel{shape}{\longrightarrow} sh \{ W,Z\}$$
\noindent {\it is a bijection.} \\

Fact 2 follows from
the proof of Lemma 3.1 in \cite{Dydak2} (the given proof is for homology but the same proof works 
for homotopy) and the Whitehead Theorem in pro-homotopy
(see \cite{Dydak1}).\\
%A continuous map $f:X\ra Z$, from a cellular disc $X$ to a space $Z$, is {\it admissible} if, for all $u\in \D^k$, we have that
%$f|_{X_u}:X_u\ra Z$ is an embedding, and $f(X_u)$ has open neighborhoods $U\sbs V$ in $Z$ such that $U$ deformation retracts, inside $V$,  to $f(X_u)$.
%Note that we do not demand the deformation retracts to depend continuously on $u$.\\

\noindent {\bf Proposition 1.1.} {\it Let $(X,\bS^{k-1})$ be
cellular $k$-disc pair with $X/\bS^{k-1}$ finite dimensional. Let $f:(X,\bS^{k-1})\ra (Z, z_0)$, $z_0\in Z$,  where $Z$ is  $LC^\infty$.
%a uniformly weakly locally contractible metric space.
%be an admissible map.
If $\pi_k(Z,z_0)=0$, $k\geq 1$,  then $f$ is null-homotopic rel $\bS^{k-1}$.}\\

\noindent {\bf Proof.}  Let $\eta:X\ra\D^k$ be the map that defines the cellular disc $X$. Write $W=X/\bS^{k-1}$.
The map $\eta$ induces a map $\eta': W\ra \D^k/\bS^{k-1}=\bS^k$.
The map $\eta'$ is a cell-like map because $\p X/\bS^{k-1}$ is homeomorphic to
$\bS^{k-1}\times X_0/\bS^{k-1}\times\{x_0\}$, hence contractible. Moreover, by hypothesis, the space $W$ is finite dimensional. Therefore, by fact 1 above,  the map $\eta'$ is a
shape equivalence, that is, an equivalence in the shape category.
%Also, the map $f$ factors through  a map $f':W\ra Z$, and since $\eta'$ is a shape equivalence there is a shape morphism
%$g\in\sh\{ \bS^k ,Z\}$ such that $g\eta'=f'$.
Consider the following commutative diagram:

$$
\begin{array}{ccc}
[ \bS^k,Z]&\stackrel{shape}{\longrightarrow}&sh\{ \bS^k,Z\}\\ \\
 (\eta')^*\downarrow \,\,\,\,\,\,\,\,\,\,\,\,\,\,\,\,& & \,\,\,\,\,\,\,\,\,\,\,\,\,\,\,\,\,\downarrow   (\eta')^* \\ \\

[  W ,Z]&\stackrel{shape}{\longrightarrow}& sh\{ W ,Z\}
\end{array} 
$$

\noindent where $(\eta')^*$ is induced by composition with $\eta'$. By fact 2 above both horizontal arrows are bijections. Also, since
$\eta'$ is a shape equivalence the right vertical arrow is also a bijection. Hence the left vertical arrow is also a bijection.
But $[\bS^k, Z]=\pi_k(Z,z_0)=0$, therefore $[W,Z]$ consists of a single element. This proves the Proposition.\\

\noindent {\bf Proposition 1.2.} {\it Let $Z$ be $LC^\infty$
%a uniformly weakly locally contractible space 
and $f:\bS^{k-1}\ra Z$. If $f$ extends to a 
cellular $k$-disc pair $(X,\bS^{k-1})$, then $f$ extends to $\D^k$.}\\

%The proof is given in Appendix B.\\

\noindent {\bf Proof.}
 Let $\eta:X\ra\D^k$ be the map that defines the cellular disc $X$.  The map $\eta$ is a cell-like map hence it induces a isomorphisms of all $i$-th
homotopy pro-groups.
Therefore all of these pro-groups are trivial for $i>0$. It follows that the inclusion $\iota:\bS^{k-1}\ra X$ represents zero in the $(k-1)$ homotopy pro-group of $X$. Consequently $f\iota$ represents zero in the
$(k-1)$ homotopy pro-group of $Z$. By fact 2 above $f\iota$ represents zero in $\pi_{k-1}(Z,z_0)$.  This proves the Proposition.\\

\noindent {\bf Proposition 1.3.} {\it  Every principal (locally trivial) $\bS^1$-bundle over a 
finite dimensional cellular disc is trivial.}\\

\noindent {\bf Proof.} Such bundles are in one-to-one correspondance with 
$[X,\C P^\infty]$, where $X$ is the cellular disc base space. Consider the following commutative diagram

$$
\begin{array}{ccc}
[ X,\C P^\infty ]&\ra& sh\{ X,\C P^\infty  \}\\
\uparrow&&\uparrow\\  

[  \D^k,\C P^\infty ]&\ra& sh\{ \D^k, \C P^\infty \}
\end{array}$$

The two horizontal maps are bijections because of Fact 2, and the right hand vertical map is also a bijection because of Fact 1. Since 
$[\D^k,\C P^\infty ]$ consists of a single point, so does 
$[X,\C P^\infty ]$. This proves the Proposition.\\

\vspace{.6in}

\noindent {\bf \large B. $C^k$-Convergence of $g$-Geodesics, with Varying $g$.}\\

Consider Riemannian metrics $g$ on a fixed manifold.
We need to study how $g$-geodesics behave when the Riemannian metric $g$ changes. We are interested in their $C^k$-convergence.
In this section $U$ denotes an open set of $\R^n$.\\

\noindent {\bf Proposition 1.4.} {\it  Let $S=\{ g^a=(g_{ij}^a)\}_{a\in A}$ be a collection of Riemannian metrics on $U$.
Let ${\bf X}=\{ {\bf x}\, /\, {\bf x}$ is a unit speed $g^a$-geodesic$, \,a\in A\}$.
Assume that the set $\{ det\, g^a(x)\,/\, x\in U,\,a\in A\}$ is bounded away from zero. Then if $S$ is $C^k$- bounded for some finite $k\geq 0$, then
the set of all derivatives $\frac{d^l x_i}{dt^l}(t)$, $1\leq l\leq k+1$, ${\bf x}=(x_1,...x_n)\in {\bf X}$, 
$t\in$ Domain of {\bf x}, is bounded.}\\

\noindent {\bf Remarks.}

\noindent {\bf 1.} The Riemannian metrics in  $S$ are not assumed to be complete.

\noindent {\bf 2.} The geodesics in ${\bf X}$ are defined on any interval.

\noindent {\bf 3.} Here ``$S$ is $C^k$-bounded" means
that for $0\leq l\leq k$, all $l$-partial derivatives of the $g_{ij}^a$ are bounded.

\noindent {\bf 4.} Note that the conclusion of the Lemma is weaker than  ``${\bf X}$ is $C^{k+1}$-bounded" (which imples $C^0$-boundedness).
Indeed, if the open set $U$ is not bounded, then ${\bf X}$ is not $C^0$-bounded, hence not $C^k$-bounded either.\\

\noindent {\bf Proof.} We denote by $(g^{ij}_a)$ the matrix inverse of $g^a=(g^a_{ij})$. First note that, since
$\{ det\, g^a(x)\,/\, x\in U,\,a\in A\}$ is bounded away from zero and $S$ is $C^k$- bounded, we have that
all $l$-partial derivatives of the $g^{ij}_a$, $0\leq l\leq k$,  are bounded.
Moreover the set $\{ \, |v|\,:\, g^a(x)(v,v)=1,\, 
v\in\R^n,\,x\in U,\,\, a\in A  \}$
is bounded. (Here $|v|$ is the Euclidean length $\langle v, v\rangle^{1/2}$.) Hence, the set of Euclidean lengths of the velocity vectors of
unit speed geodesics is bounded. Therefore, the set $$\bigg\{\, \frac{d\, x_i(t)}{dt}, \,{\bf x}(t)=(x_1(t),...,x_n(t)) {\mbox{ is a unit speed }}
g^a{\mbox{-geodesic}},\, a\in A,\, t\in {\mbox{Domain of }}{\bf x}\,\bigg\}$$
\noindent is bounded. This proves the Proposition for $k=0$.\\

Assume $S$ is $C^1$-bounded. Let ${\bf x}(t)=(x_1(t),...,x_n(t))$ be a unit speed $g^a$-geodesic. 
Then the $x_i$'s satisfy a second order ODE of the form
$$\frac{d^2 x_i}{dt^2}=\Phi \bigg(\, \frac{dx_j}{dt},\, \Gamma_{st}^r({\bf x})\,\bigg) $$
\noindent where $\Gamma_{st}^r=(\Gamma_{st}^r)^a$ are the Christoffel symbols of the metric $g^a$ and the function $\Phi$ is a polynomial
function independent of ${\bf x}$ and $a\in A$. But the Christoffel symbols $(\Gamma_{st}^r)^a$ can be written canonically as a polynomial expression on
the $g^{ij}_a$ and the first partial derivatives of the $g_{ij}^a$. Since all these terms are bounded we conclude that the set of all
Christoffell symbols $(\Gamma_{st}^r)^a$ is bounded. Therefore the set 
$$\bigg\{\, \frac{d^2\, x_i(t)}{dt^2}, \,{\bf x}(t)=(x_1(t),...,x_n(t)) {\mbox{ is a unit speed }}
g^a{\mbox{-geodesic}},\, a\in A,\, t\in {\mbox{Domain of }}{\bf x}\,\bigg\}$$
\noindent is also bounded. This proves the Proposition for $k=1$. \\

Assume $S$ is $C^2$-bounded. We differentiate the geodesic equation above to obtain the
third order ODE
$$\frac{d^3 x_i}{dt^3}=\Psi \bigg(\, \frac{dx_j}{dt}, \, \frac{dx^2_j}{dt^2},\,\Gamma_{st}^r({\bf x}), \,\frac{\p^k}{\p x^k}\Gamma_{st}^r({\bf x})\,\bigg) $$
\noindent which is satisfied by any  ${\bf x}=(x_1,...x_n)\in {\bf X}$. Since $\Psi$ is a universal polynomial, and $\Psi$ is applied
to a set of bounded variables we conclude that the Proposition holds for $k=2$. Proceeding in this way we prove the Proposition for any $k\geq 0$.
This proves the Proposition.\\

In the next two Propositions we use the following notation. For a Riemannian metric $g_0$ and sequence of Riemannian metrics $\{g_n\}$ on $U$ 
we write $g_n\stackrel{C^k}{\ra} g_0$ to express uniform $C^k$-convergence
on compact supports. Also, for $p\in U$, $v\in\R^n$ we denote by $\alpha(p,v,g)$ the $g$-geodesic
with value $p$ at zero, and velocity $v$ at zero.
Also $\alpha(p_n,v_n,g_n)\stackrel{C^{k}}{\ra}\alpha(p,v,g_0)$
means convergence on any
closed interval $[a,b]$ where all paths are defined. (Note that in this case there is $\epsilon >0$ such that $\alpha(p,v,g_0)$ and all 
$\alpha(p_n,v_n,g_n)$ are defined on $[-\epsilon, \epsilon]$.)\\

\noindent {\bf  Lemma 1.5.} {\it If  $g_n\stackrel{C^1}{\ra} g_0$, $p_n\ra p$, $v_n\ra v$, then 
$\alpha(p_n,v_n,g_n)\stackrel{C^{1}}{\ra}\alpha(p,v,g_0)$ }\\

\noindent {\bf Proof.} $C^1$-Convergence follows from the general theory of first order ODE with parameters. This proves the Lemma.\\

\noindent {\bf Proposition 1.6.} {\it  Let  $g_n\stackrel{C^k}{\ra} g_0$, $k\geq 1$,  and $\alpha_n(t)$, $t\in[a,b]$, be $g_n$-geodesics such that
$\alpha_n\stackrel{C^0}{\ra}\alpha_0$. Then
$\alpha_n\stackrel{C^{k+1}}{\ra}\alpha_0$. }\\

\noindent {\bf Proof.} It is enough to prove $C^1$-convergence because then the $C^k$-convergence, $k\geq 2$, 
follows using the same argument used in the proof of Proposition 1.4 involving the $\Phi$, $\Psi$,...functions.
But if $\alpha_n$ does not $C^1$-converge to $\alpha_0$ we arrive, using Lemma 1.5, to a contradiction. This proves the Proposition.
\vspace{.6in}

\noindent {\bf \large C.  Sets of Parallel Lines in a Hadamard Manifold.}\\

Let $H=H^n$ be a Hadamard manifold and $\cL$ a set of parallel geodesic lines in $H$. We assume that $\cL$ is {\it ribbon convex}, i.e. if $\ell_0,\,\ell_1\in\cL$ then $\ell\in\cL$, for every $\ell$ contained in the flat ribbon bounded by $\ell_0$ and $\ell_1$ (for the existence of the flat ribbon see \cite{LawYa}). 
Write $L=\bigcup\cL$. The Flat Ribbon Theorem of A. Wolf \cite{LawYa} implies that
$L$ is a convex set.  We choose one of the two points at infinity determined by any $\ell\in\cL$. This choice ``orients" all lines $\ell\in\cL$ and we can now make   
the real line $\R$ act isometrically on $L$ by translations: for $t\in\R$ and $p\in\ell\in\cL$, $t.p=q$, where $q\in\ell$, and $q$ is obtained from
$p$ by a $t$-translation.\\

Now, fix  $p\in \ell_0 \sbs L$. Let $\ell\in \cL$. Since $\ell_0$ and $\ell$ bound a flat ribbon there is a unique point $p_\ell\in\ell$ which
is the closest to $\ell_0$ and the geodesic segment $[p,p_\ell]$ is perpendicular to both $\ell_0$ and $\ell$. Write $K=\{\, p_\ell\,\,   |\,\, \ell\in\cL\,\}$.
Note that $K\cap \ell=p_\ell$.\\

\noindent {\bf Proposition 1.7.} {\it  The set $K$ is convex.}\\

%We need the following Lemma.\\

%\noindent {\bf Lemma 1.8} {\it Let $\ell_1, \ell_2,  \ell_3 \in\cL$ and $p_1, p_4\in\ell_1$,  $p_2\in\ell_2$, %$p_3\in\ell_3$ such that
% $[p_1,p_2]$ is perpendicular to $\ell_1$, $\ell_2$, $[p_2,p_3]$ is perpendicular to $\ell_2$, $\ell_3$, 
%$[p_3,p_4]$ is perpendicular to $\ell_3$, $\ell_1$. Then $p_1=p_4$.}\\ 

%\noindent {\bf Proof of the Lemma.} We can assume all lines $\ell_i$ to be different. Let $S$ be the union of the three ribbons bounded by each pair $\ell_i$,
%$\ell_j$. Then $S$, with the intrinsic metric (i.e the distance is the smallest length of a path joining two points), is isometric to the cylinder
%$\bS^1(r)\times \R$, where $\bS^1(r)\sbs \R^2$ is the circle of radius $r$. But if $p_1\neq p_4$, the $\ell_i$ do not globally minimize length in $S$,
%hence they do not minimize length in $H$, which is a contradiction. This proves the Lemma.\\

%\noindent {\bf Proof of Proposition 1.7.} Let $p_1, p_2\in K$, $p_i\in\ell_i$. We prove $[p_1,p_2]\sbs K$.
%By Lemma 1.5 $[p_1,p_2]$ is perpendicular to $\ell_1$, $\ell_2$. Let $p_3\in[p_1,p_2]$, $p_3\in\ell_3$. Note that 
%$[p_1,p_3]$ is perpendicular to $\ell_1$, $\ell_3$. By Lemma 1.8
%$[p,p_3]$ is perpendicular to $\ell_0$, $\ell_3$, hence $p_3=p_{\ell_3}$. Therefore $p_3\in K$. This proves the Proposition.\\

This lemma is proved in \cite{FJ}.\\

Consider the map $K\times\R\ra L$, $(p, t)\mapsto t.p_\ell$. Since $K$ is convex this map is an isometry, where we consider
$K\times\R$ with the metric $d( (p,t), (p',t'))=\sqrt{d_H(p, p')^2+d_{\R} (t,t')^2}$. The inverse of this map is the map $(\pi, T)$,
where $\pi (x)=p_\ell$, $x\in \ell$, is the projection onto $K$, and $T(x)$ is the (oriented) distance between $x$ and $p_\ell$.\\

\noindent {\bf Corollary 1.8.} {\it Assume $L$ is a closed subset of $H$
and that $K$ is compact. Then 
$K$ is homeomorphic to the closed $k$-disc with smoothly (locally) totally geodesic embedded interior.}\\

\noindent {\bf Proof.} Proposition 1.7 and
Theorem 1.6 of \cite{Cheeger-Gromoll}, p. 418, imply that
$K$ is homeomorphic to a compact, contractible $k$-manifold, $0\leq k\leq n-1$. Moreover the inclusion $K\hookrightarrow H$ restricted to the (manifold) interior of $K$ is smooth and 
(locally) totally geodesic. This proves the corollary.
\vspace{.6in}

\noindent {\bf \large D.  Sets of Homotopic Closed Geodesics.}\\

Let $Q=\bS^1\times\R^{n-1}$, with a complete nonpositively curved Riemannian metric $g$.
Write $\iota :\bS^1=\bS^1\times \{ 0\}\hookrightarrow Q$ for the inclusion and 

$$\Omega=\{\,\,\alpha\in C^\infty (\bS^1, Q)\,\, / \,\, \alpha \simeq \iota\,\, \}$$

\noindent with the $C^k$ topology, $0\leq k\leq\infty$.
Note that $\bS^1$ acts freely on $\Omega$ by $z.\alpha (w)=\alpha (zw)$, for $z, w\in\bS^1\sbs\C$.
Write $\Sigma=\Omega/\bS^1$. It is straightforward to verify that the quotient map $\Omega\ra\Sigma$ is a 
(locally trivial) principal $\bS^1$-bundle. Note that $\R$ also acts on $\Omega$ by $z=x.\alpha (w)=\alpha (e^{2\pi i x}w)$, for $x\in\R$, $w\in\bS^1$.
Moreover, we also get $\Omega/\R=\Sigma$.
Let $\cC=\cC_g$ be the set of all parametrized closed geodesics
homotopic to the inclusion, i.e. 

$$\cC=\{\,\,\alpha\in \Omega\,\,/\,\,\alpha \,\, {\mbox{is a}}\,\,\, g-{\mbox{geodesic}}\,\, \}$$

\noindent Note that every $\alpha\in \cC$ is an embedding. Moreover, by the Flat Ribbon Theorem
of J. A. Wolf \cite{LawYa} any two elements in $\cC$ either have the same image (hence lie in the same $\bS^1$-orbit) or have disjoint images.\\

Let $\cA=\cA_g$ be the image of $\cC$ by the bundle map $\Omega\ra\Sigma$. That is, $\cA$ 
is the set of all unparametrized closed geodesics homotopic to the  inclusion.
Assume that\\

{\bf (a)} $\cC$ is non-empty.

{\bf (b)} $\cC$ is $C^0$-bounded. Equivalently, the set $C=\bigcup \cC$
is contained in a compact set.\\

\noindent {\bf Proposition 1.9.} {\it Under these assumptions $\cA$ is homeomorphic to a closed $k$-disc, $k\leq n-1$.}\\

\noindent {\bf Remark.} All topologies $C^k$, $0\leq k\leq \infty$, induce the same topology on
$\cC$ and $\cA$.\\

\noindent {\bf Proof.} Let $H$ be the universal cover of $Q$. Then $H$ is a Hadamard manifold, the infinite cyclic group $\Z$ acts freely 
by isometries on $H$ and $Q=H/\Z$. 
%Fix $\alpha_0\in\cC$. 
%Let $\beta_0:\R\ra H$ be a lifting of $\alpha_0:\bS^1\ra Q$. Write $\ell_0=\beta_0(\R)$. 
Let $\cL$ 
be the set of all lines in $H$ which cover elements in $\cA$; i.e. all lifts to $H$ of unparametrized
closed geodesics homotopic to $\iota$. It is straightforward to check that $\cL$ is ribbon-convex and $L=\bigcup\cL$
is closed. Let $K$ be constructed from $\cL$ as in section C. Using projections we can construct,
in the obvious way, a one-to-one, onto $\cA$, continuous map $K\ra \cA$. This proves the Proposition.\\

\noindent {\bf Proposition 1.10.} {\it The space $\Omega$ deformation retracts to $\cC$.}\\

\noindent {\bf Proof.} 
Let $\Omega_\iota Q$ denote the space of all based loops which are based homotopic to $\iota$. Then
we have a fibration $\Omega_\iota Q\ra\Omega\ra Q$, where the last map is
the evaluation map (at, say, $1\in\bS^1$). Since $\Omega_\iota Q$ is
contractible and  $Q\sim\bS^1$ it follows that $\cC\hookrightarrow\Omega$ is a homotopy equivalence. This proves the proposition.\\

Since $\cA$ is a disc, the bundle $\Omega\ra\Sigma$ restricted to $\cA\sbs \Sigma$ is trivial.
Hence $\cC$ is homeomorphic to $\cA\times \bS^1$. Let $s:\cA\ra\cC$ be any section of this bundle
(equivalently, a lifting of the identity $1_\cA$). Write $\cB=s(\cA)\sbs\cC\sbs\Omega$.\\

\noindent {\bf Proposition 1.11.} {\it 
Let $V$ be an open neighborhood of $\cB$ in $\Omega$. There the is an open neighborhood $U\sbs V$ of $\cB$
in $\Omega$ such that $U$ deformation retracts, inside $V$, to $\cB$}.\\

\noindent {\bf Proof.} Let $h_t$, $h_0=1_\Omega$, $h_1:\Omega\ra\cC$, be a deformation retract. 
Since $\cC$ is homeomorphic to $\cA\times \bS^1$, we can find an open neighborhood $W$ of $\cB$ in $\cC$ that deformation retracts to $\cB$.
And we can assume that this deformation retract happens inside $V$. Denote this deformation retract by $f_t$.
Let $U$ be an open neighborhood of $\cB$ in $\Omega$ small enough so that
$h_t(U)\sbs V$ and $h_1(U)\sbs W$. Then our desired deformation retract is the concatenation of the $h_t$ with the $f_t$.
This proves the Proposition.\\

We will need the following Lemma in the next section.\\

\noindent {\bf Lemma 1.12.} {\it Assume the metric $g$ on $Q$ satisfies assumptions (a) and (b) above. 
Let $\pounds$ be the length of a (hence all) closed $g$-geodesic homotopic to $\iota$. Then there is a bounded set
$R\sbs Q$ such that if the image of an $\alpha\in \Omega$ is not contained in $R$ then the $g$-length of $\alpha$ is
larger than $1+\pounds$.}\\

\noindent {\bf Proof.} Suppose not. Then there is a sequence $\alpha_n$ in $\Omega$, such that: 
(1) $x_n=\alpha_n(1)$ goes to infinity and (2) all $\alpha_n$ have length $\leq 1+\pounds$. 
Fix $\alpha_0\in\cC$ and write $x=\alpha_0(1)$. Let $s_n$ be a geodesic segment $[x,x_n]$ such that $d(x, x_n)$ is its length
and write $s_n(t)=exp_x(tv_n)$, for some
unit length vector $v_n\in T_xQ$. We can assume $v_n\ra v$, where $v$ also has unit length. Write $s(t)=exp_x(tv)$. Let $H$ be the 
universal cover of $Q$. Fix a lift $\beta_0:\R\ra H$ of $\alpha_0$. Write $y=\beta_0(0)$ and $z=\beta_0(1)$. Let $s_n', s'$ be liftings of
$s_n$ and $s$ beginning at $y$ and  $s_n'', s''$ be liftings of
$s_n$ and $s$ beginning at $z$, respectively. Note that the endpoints of $s_n'$ and $s_n''$ can be joined by a lifting of $\alpha_n$, hence
their distance lies in the interval $[\pounds, 1+\pounds]$. Therefore $d_H(s_n'(t), s_n''(t))\in [\pounds, 1+\pounds]$.
It follows that $d_H(s'(t), s''(t))\in [\pounds, 1+\pounds]$ for all $t\geq 0$. But the function $t\mapsto d_H(s'(t), s''(t))$ is convex
with minimum value at $t=0$, thus it cannot be a bounded function unless it is constant. But this contradicts assumption (b). This proves the Lemma.

\vspace{.6in}

\noindent {\bf \large E.  Sets of Homotopic Closed $g$-Geodesics, with Varying $g$.}\\

Let $Q, \Omega, \Sigma$ be as in section D. We denote by $\met^{sec\leq 0}(Q)$ the space of all complete nonpositively curved Riemannian metrics
on $Q$, with the weak smooth topology (i.e the union of the weak $C^s$ topologies, which are the topologies of the $C^s$-convergence on compact sets).
Let $\sigma:\D^k\ra\met^{sec\leq 0}(Q)$ be continuous. Write $g_u=\sigma(u)$. Using the methods and the notation of section
D, for each $g_u$ we obtain $\cA_u$, $\cB_u$, $\cC_u$.
Write $C_u=\bigcup\cC_u$.
In what follows we assume that all $g_u$ satisfy assumptions (a) and (b) of section D. 
In particular, for each $u$ we get a positive number $\pounds (u)$ which is
the length of an element in $\cC_u$, that is, the length of a $g_u$-geodesic homotopic to the inclusion $\bS^1\ra Q$.\\

% \noindent {\br Lemma 1.13.} {\it The set $\{\pounds (u)\}_{u\in\D^k}$
%is bounded.}\\

%\noindent {\bf Proof.} Suppose not. Then there is $u_n\ra u$ with
%$\pounds_n=\pounds_{u_n}\ra infty$. Let $\alpha$ a $g_u$-geodesic
%with $g_u$-length $\ell_{g_u}(\alpha)=\pounds_u$. For $n$ large
%$\ell_{g_u}(\alpha)$ and $\ell_{g_{u_n}}(\alpha)$ are close hence we get
%$\ell_{g_{u_n}}(\alpha)<\pounds_{u_n}$, for $n$ large. This contradicts
%the minimality of the $g_{u_n}$-length of $g_{u_n}$-geodesics.
%This proves the lemma.\\

\noindent {\bf Lemma 1.13.} {\it The map $\pounds:\D^k\ra (0,\infty)$ is upper semi-continuous.}\\

\noindent {\bf Proof.} This follows from the following facts: (1) $\pounds (u)$ is the smallest possible length of a curve
homotopic to the inclusion $\bS^1\ra Q$, and 
(2) for any closed curve $\alpha$, the $g_u$-length of $\alpha$ is close to
the $g_v$-length of $\alpha$, provided $u$ is close to $v$. This proves the Lemma.\\

\noindent {\bf Lemma 1.14.} {\it Let $u_n\ra u$ in $\D^k$.
Then there is a compact set $S$ of $Q$ and a sequence $\alpha_n\in\cC_{n}$ such that $\alpha_n\sbs S$, for $n$ sufficiently
large.}\\

\noindent {\bf Proof.} Let $R$ be as in lemma 1.12 for $g=g_u$
and assume that $R$ is closed. Let $S$ be any compact of $Q$ with 
$R\sbs int\, S$. We claim that there is a sequence $\alpha_n\in\cC_{n}$ such that $\alpha_n\sbs S$, for $n$ sufficiently large. This would prove the lemma.
Suppose not. Then we can assume, by passing to a subsequence,  that for every $\alpha_n\in \cC_n$ we have $\alpha_n\not\subset S$.
Write $\pounds=\pounds(u)$, $g_n=g_{u_n}$ and let $\alpha\in \cC_u$. Then the $g$-length $\ell_g(\alpha)$ of
$\alpha$ is $\pounds$. Therefore $\ell_{g_n}(\alpha)$ is close to $\pounds$. For each
$n$ let $\alpha^t_n$ be a homotopy with $\alpha_n^0=\alpha$, $\alpha_n^1\in\cC_{u_n}$
and $\ell_{g_n}(\alpha_n^t)\leq\ell_{g_n}(\alpha_n^s)$, for $t>s$. That is the deformation
$t\mapsto \alpha_n^t$ begins in $\alpha$, ends in a $g_n$-geodesic, and is $g_n$-length
non-increasing. (Such a deformation can be done in the usual way using evolution
equations or using a polygonal deformation.) 
Note that, by hypothesis, $\alpha_n^1\not\subset S$ and $\alpha=\alpha_n^0\sbs R\sbs S$. This together with the continuity
of the deformation imply that there is $s=s_n$ such that $\beta_n=\alpha_n^s\not\sbs R$ and
$\beta_n\sbs S$. But lemma 1.12 together with the convergence $g_n\ra g$ imply that
$\ell_{g_n}(\beta_n)>1/2+\pounds$ when $n$ is sufficiently large (see remark below). This contradicts the fact that
the deformation $t\mapsto \alpha_n^t$  is $g_n$-length
non-increasing. This proves the lemma.\\

\noindent {\bf Remark.} In the proof above we are using the following fact:
if $|g_n-g|_g\leq\delta$ then for any $PD$ path $\alpha$ we have $|\ell_{g_n}(\alpha)-\ell_{g}(\alpha)|\leq
\frac{\delta}{1-\delta}\ell_{g_n}(\alpha)$. This fact follows from the definition of length (using integrals) and
the triangular inequality.\\

\noindent {\bf Corollary 1.15.} {\it The map $\pounds:\D^k\ra (0,\infty)$ is continuous.}\\

\noindent {\bf Proof.} Let $u_n\ra u$ and $S$ be as in 1.14. Hence there are $\alpha_n\in \cC_{u_n}$
with $\alpha_n\sbs S$. Since $g_{u_n}\ra g_u$ uniformly on $S$
we get from lemma 1.13 (and the remark above) that $\pounds(u_n)=\ell_{g_{_{u_n}}}(\alpha_n)$ is close to $\ell_{g_u}(\alpha_n)\geq \pounds (u)$. This shows $\pounds$ is lower semi-continuous.
This proves the corollary.\\

\noindent {\bf Proposition 1.16.} {\it The set $\bigcup_{u\in\D^k}\cC_u$ is $C^0$-bounded.}\\

That is, the set of all $g_u$-geodesics lie at bounded distance from, say, the inclusion $\iota$, for all $u\in\D^k$.\\ 

\noindent {\bf Proof.} Suppose not. Then there are $u_n\ra u$ and $\alpha'_n\in\cC_n=\cC_{u_n}$
with $\alpha'_n$ going to infinity, i.e. $\alpha'_n\not\sbs K$, for any given compact $K$, provided $n$ is large.
Let $S$ be as in 1.14. Hence there are $\alpha_n\in \cC_{u_n}$
with $\alpha_n\sbs S$. Let $S'$ be a compact such that $S\sbs int\, S'$. Since $\alpha_n$ and $\alpha'_n$
bound a flat two dimensional cylinder (in the $g_n=g_{u_n}$ metric) we can find 
$\beta_n\in\cC_n$ with $\beta_n\in S'$ and $\beta_n\not\sbs S$. By corollary 1.15 we can assume
$\pounds(u_n)\leq 1/2+\pounds(u)$. On the other hand, by 1.12 and the uniform convergence
$g_n\ra g_u$ on $S'$ we have $\pounds (u_n)=\ell_{g_n}(\beta_n)$ is close to $1+\pounds(u)$.
This is a contradiction. This proves the proposition.\\

\noindent {\bf Proposition 1.17.} {\it The set $\bigcup_{u\in\D^k}\cC_u$ is $C^k$-bounded, for any $k$,
$0\leq k<\infty$.}\\

\noindent {\bf Proof.} The Proposition follows from Proposition 1.4 by considering $Q$ as an open set of $\R^n$. Note that, by 1.16, we can work on an open set with compact closure, hence all
required quantities will be bounded. Note also that in 1.4 the geodesics are assumed to have
speed one, but the geodesics in $\cC_u$ have speed $\pounds (u)/2\pi$. This  can be fixed by a rescaling
of geodesics and using the fact that (by 1.15) the set $\{\pounds (u)\}_{u\in\D^k}$ is bounded and bounded away from zero. This proves the proposition.\\

\noindent {\bf Proposition 1.18.} {\it Let $V$ be an open neighborhood of $\cC_u$ in $\Omega$. Then, for $u'\in \D^k$ close enough to $u$,
$\cC_{u'}\sbs V$.}\\

\noindent {\bf Proof.} Suppose not. Then there are $u_n\ra u$ and $\alpha_n\in\cC_{u_n}$ with $\alpha_n\notin V$.
Since $u_n\ra u$, Proposition 1.17 (recall we are using the weak smooth topology) implies that the set $\{\alpha_n\}$ is $C^0$-equicontinuous,
where we consider $Q$ with metric $g_u$. Moreover we can assume all $\alpha_n$ to be Lipschitz with the same constant.
Proposition 1.17 also says that the set $\{ \alpha_n\}$ is bounded. By Arzela-Ascoli Theorem we can assume the $\alpha_n$ $C^0$-converge 
to a Lipschitz $\alpha:\bS^1\ra Q$. Since $\alpha$ is Lipschitz its length is finite (the length defined as
$sup \sum d_{g_u}(\alpha (z_i), \alpha (z_{i+1}))$, the $sup$ taken over all partitions of $\bS^1$). Moreover it is straightforward to
show that $g_{u_n}$-lengths of the $\alpha_n$ converge to the $g_u$-length of
$\alpha$. Proposition 1.15 implies now that $\alpha$ has minimal $g_u$-length, hence it is smooth and $\alpha\in\cC_u$. Now using proposition 1.6 we see that  $\alpha_n\in V$ for $n$ large enough.
This contradiction proves the proposition.\\

Define $Y=\coprod_{u\in\D^k} \{u\}\times \cA_u\sbs \D^k\times \Sigma$, that is
$ Y=\{\,\, (u, a)\,\, |\,\, u\in\D^k,\,\,\, a\in \cA_u\,\,   \}$.
Define also $Z=\coprod_{u\in\D^k} \{u\}\times \cC_u\sbs \D^k\times \Omega$. Each $C^k$-topology on  $\Omega$, $0\leq k\leq \infty$,
induces a $C^k$-topology on $Z$.\\

\noindent {\bf Proposition 1.19.} {\it All   \,$C^k$-topologies on $Z$ coincide.}\\

\noindent {\bf Proof.} This follows from Proposition 1.6. This proves the Proposition.\\

\noindent {\bf Proposition 1.20.} {\it The space $Z$ is compact and metrizable.}\\

\noindent {\bf Proof.} The space $Z$ is certainly metrizable. Let $\{(u_n,\alpha_n)\}$ be a sequence in $Z$.
We can assume $u_n\ra u$.
By Proposition 1.17, the sequence $\{ \alpha_n\}$ is $C^k$-bounded, for all $k\geq 0$. In particular it is
$C^1$-bounded. Therefore the sequence $\{ \alpha_n\}$ is equicontinuous.
Moreover we can assume all $\alpha_n$ to be Lipschitz with the same constant.
Proposition 1.16 says that the set $\{ \alpha_n\}$ is $C^0$-bounded. By Arzela-Ascoli Theorem we can assume that $\{\alpha_n\}$ $C^0$-converges 
to a Lipschitz $\alpha\in\Omega$. Since $\alpha$ is Lipschitz its length is finite, where the length is defined as
$sup \sum d_{g_u}(\alpha (z_i), \alpha (z_{i+1}))$, the $sup$ taken over all partitions of $\bS^1$. Using
this definition of length it is straightforward to
show that $g_{u_n}$-lengths of the $\alpha_n$ converge to the $g_u$-length of
$\alpha$. Corollary 1.15 implies now that $\alpha$ has minimal $g_u$-length, hence it is smooth and $\alpha\in\cC_u$.
Therefore $\{(u_n,\alpha_n)\}$ converges to $(u,\alpha)\in Z$. This proves the Proposition.\\

It follows from Proposition 1.6 that $Y$ is compact. The space $Y$ is also Hausdorff.
% and second countable (if $\{B_i\}$ is a countable basis for $Z$
%the sets $\bS^1B_i$ determine a countable basis for $Y$.) 
Therefore the projection $Y\ra\D^k$ is a cellular $k$-disc (choose any base point). \\

\noindent {\bf Proposition 1.21.} {\it The space $Y$ is finite
dimensional.}\\

\noindent {\bf Proof.}  Since $Y$ is Hausdorff compact and $Z\ra Y$ is locally
trivial, the proposition follows from the following claim.\\

\noindent {\bf Claim.} {\it
Let $U\sbs Y$ be compact and such that $Z\ra Y$ is trivial over $U$. Then
$U$ is homeomorphic to a compact subset of $\R^{n+k}$. Hence $U$ is finite dimensional.}\\

\noindent{\bf Proof of the claim.} Let $U'\sbs Z$ be the image of a section of $Z\ra Y$
over $U$. Hence the restriction $U'\ra U$  is a homeomorphism.
As before we are considering $Q$ as an open set of $\R^n$. Now, just define $h:U'\ra\R^{k+n}$ as $h(u,\alpha)=(u,\alpha (1))$.
(Recall $1\in\bS^1\sbs\C$ and $\alpha:\bS^1\ra Q$.) This is a one-to-one continuous map with compact domain between metric spaces.
Hence it is a homeomorphism onto its image. This proves the claim
and Proposition 1.21.\\

By Propositions
1.3 and 1.21 the $S^1$-bundle $Z\ra Y$ is trivial (thus $Z$ is homeomorphic to $Y\times\bS^1$). Take a section $Y\ra Z$ of this bundle and let
$X$ be the image of $Y$ by this section. Write $\eta : X\ra \D^k$ for the projection.
Then $\eta:X\ra \D^k$ is a cellular $k$-disc. Note that $X$ is formed by honest {\it parametrized} $g_u$-geodesics, not unparametrized ones, like
the ones in $Y$. And since $X$ and $Y$ are homeomorphic, Proposition 1.21 has the following corollary.\\

\noindent {\bf Corollary 1.22.} {\it The space $X$ is a finite dimensional space.}\\

Recall that the cellular discs $Y$ and $X$ were constructed from a map $\sigma:\D^k\ra\moo(M)$.
Now assume that $\sigma|_{\bS^{k-1}}$ is {\it constant mod $DIFF(M)$}; that is
$\sigma|_{\bS^{k-1}}$ factors through the map $\Lambda_{g_0}$ in \cite{FO6}.
(This map is just the orbit map
$\Lambda_{g_0}:DIFF(M)\ra \moo(M)
$ given by $\phi\mapsto \phi_*g_0$.)
Hence for $u\in\bS^{k-1}$ we can write $\sigma_u=(\phi_u)_*g_0$, for
some continuous $u\mapsto\phi_u$, and it follows that
$\cA_u=\phi_u(\cA_{u_0})$, $u\in\bS^{k-1}$.
Therefore we can write $\p Y=\bS^{k-1}
\times\cA_{u_0}$ and we can consider $\bS^{k-1}\sbs Y$ by choosing
any element in $\cA_{u_0}$. Analogously for $X$.
Hence we obtain cellular $k$-discs pairs $(Y,\bS^{k-1})$, $(X,\bS^{k-1})$.\\

\noindent {\bf Corollary 1.23.} {\it 
Assume $\sigma|_{\bS^{k-1}}:\bS^{k-1}\ra\moo(M)$ is constant mod $DIFF(M)$. Then the space $X/\bS^{k-1}$ is a finite dimensional space.}\\

\noindent {\bf Proof.} From the proof of 1.21 we have that the embedding $h$ induces an
embedding $X/\bS^{k-1}\ra \R^{k+n}/h(\bS^{k-1})$ which is clearly finite dimensional.\\

\vspace{.8in}

\noindent {\bf \Large  Section 2. Proof of the Main Theorem.}\\

In this section we shall use the notation and results given \cite{FO6}
to prove the nonpositively curved version of the Main Theorem in \cite{FO6}.
In turn this version follows from the nonpositively version of
Theorem 1 of \cite{FO6} together with Theorem 2 of \cite{FO6}.
Right before Theorem 1 of \cite{FO6} the following diagram is given:

$$\begin{array}{ccccc} DIFF( \, (\bS^1\times\bS^{n-2})\times I ,\p\, )&\stackrel{\Phi}{\ra}& DIFF(M)&
\stackrel{\Lambda_g}{\ra}&\mo (M)\\&&&&\\
\iota\,\, \downarrow \,\,\,\,\,& & &&\\&&&&\\  P(\bS^1\times\bS^{n-2})&&&&\end{array}$$

\noindent where $P(N)$ is the space of topological pseudoisotopies
of a manifold $N$. In \cite{FO6} it is proved that, under the relevant conditions,
$Ker \, (\, \pi_k(\Lambda_g \Phi)\,)\sbs Ker\, (\, \pi_k(\iota )\, )$.
Before we present the nonpositively curved  version of this result
we need some notation and definitions. \\
%For this we need the following variation of the
%pseudo-isotopy space. For a compact manifold $N$ let $\hP(N)$ be the space of all cellular maps
%$N\times [0,1]\ra N\times [0,1]$, that are the identity on $(N\times \{ 0\})\cup (\p N\times I)$ and
%are a homeomorphism when restricted to $N\times [0,1)$.
%We consider $\hP(N)$ with the compact-open topology. Then we have the following diagram of topological
%inclusions
%$$\begin{array} {ccccc}P(N)&\sbs & \hP(N)&&
%\\  \cap &&\cap&&  \\
%TOP(N\times [0,1])&\sbs& CELL(N\times[0,1])&\sbs& C(N\times [0,1])\end{array}$$

For a manifold $L$, $TOP(L)$ is the space of all self homeomorphisms of $L$, $CELL(L)$ is the space of cellular 
maps and $C(L)$ is the space of all continuous $L\ra L$,
with the compact-open topology.\\

\noindent {\bf Remark.} If the manifold $L$ has boundary then for a map
$f:L\ra L$ to be cellular we demand the restriction $f|_{\p L}:\p L\ra\p L$
to be cellular too. See \cite{Si}.\\

\noindent{\bf Lemma 2.1.} {\it The map $\pi_k\, P(N)\ra\pi_k \, TOP(N\times [0,1])$, $k\geq 0$,  is injective.}\\

\noindent {\bf Proof.} Let $\alpha: \bS^{k}\ra P(N)$, $\beta:\D^{k+1}\ra TOP(N\times [0,1])$, with
$\beta|_{\bS^k}=\alpha$. For $f\in TOP(N\times [0,1])$ write $f_0:N\ra N$ for its bottom, that is, for its 
restriction to $N\times \{ 0\}$. Define $\gamma:\D^{k+1}\ra TOP(N\times [0,1])$ by $\gamma (u)=(\beta(u)_0)^{-1}\times 1_{[0,1]}$.
Note that $\gamma(u)=1_{N\times[0,1]}$ for $u\in\bS^k$.
Finally define  $\beta':\D^{k+1}\ra TOP(N\times [0,1])$ by $\beta'(u)=\gamma(u) \beta(u)$. Then $\beta'|_{\bS^k}=\alpha$.
This proves the Lemma because $\beta':\D^{k+1}\ra P(N)$.\\

\noindent{\bf Lemma 2.2.} {\it Let $N$ be compact and $dim\, N\neq 3$. Then the 
 map $\pi_k \, TOP(N\times [0,1])\ra\pi_k \, CELL(N\times [0,1])$, $k\geq 0$, is an isomorphism.}\\

This is a fiber version of Siebenmann result \cite{Si}. The proof follows from
 proposition 4.1 of  B. Haver \cite{Haver} together with the fact that the closure of 
$TOP(L)$ is $CELL(L)$, $dim L\neq 4$, proved by Siebenmann \cite{Si}.
In the Lemma above (and the Corollary below), for the case $k=0$ ``isomorphism'' means ``bijection''. These two Lemmas imply:
\\

\noindent{\bf Corollary 2.3.} {\it Let $N$ be compact and $dim\, N\neq 3$. Then the 
 map $\pi_k \, P(N)\ra\pi_k \, CELL(N\times [0,1])$, $k\geq 0$, is injective.}\\

Now, assume $M$ is closed and admits a nonpositively curved metric $g$. Here is the new version of the diagram above:

$$\begin{array}{ccccc} DIFF( \, (\bS^1\times\bS^{n-2})\times I ,\p\, )&\stackrel{\Phi}{\ra}& DIFF(M)&
\stackrel{\Lambda_g}{\ra}&\moo (M)\\&&&&\\
\iota\,\, \downarrow \,\,\,\,\,& & &&\\&&&&\\  P(\bS^1\times\bS^{n-2})&&&&\\ \\
\iota'\,\, \downarrow \,\,\,\,\,& & &&\\&&&&\\ CELL(\bS^1\times\bS^{n-2}\times [0,1])&&&&
\end{array}$$

\noindent where $\iota'$ is the inclusion. 
And we will prove the following version of Theorem 1 of \cite{FO6}:\\

\noindent {\bf Theorem 2.4.} {\it Let $M$ be a closed $n$-manifold with a nonpositively
curved metric $g$. Let $\alpha$, $V$, $r$ and $\Phi=\Phi (\alpha , V, r)$ be as in \cite{FO6},
and assume that $\alpha$ in not null-homotopic. Then 
$Ker \, (\, \pi_k(\Lambda_g \Phi)\,)\sbs Ker\, (\, \pi_k(\iota'\iota )\, )$, for $k<n-5$.}\\

Note that the Corollary above says that $\pi_k(\iota')$ is injective. Therefore
we just need to prove Theorem 2.4. A few modifications have to be made
in section 2 of \cite{FO6}, since there it is assumed that the manifold $(Q,g)$ is negatively
curved. Here we just assume that it is nonpositively curved, and in addition the universal cover $\tQ$
is hyperbolic (that is, as a geodesic space $\tQ$ is $\delta$-hyperbolic, $\delta>0$).\\

\noindent {\bf Remark.}  We will use the fact that $CELL(L)$ is 
$LC^\infty$ for a compact $L$, $dim \,L\neq 3$. This follows from proposition 4.1 of \cite{Haver} and the main
results in \cite{Si} and \cite{Edwards-Kirby}.
%CHECK IN \cite{Edwards-Kirby} IF THE BOUNDARY IS ALLOWED TO MOVE.
\vspace{.6in}

\noindent {\bf \large Changes to Section 2 of \cite{FO6}.}\\
There are several facts stated for $Q$ in section 2 of \cite{FO6} that are clearly
not true for a general nonpositively curved manifold. But, since we are assuming
$\pi_1M$ word hyperbolic we can assume that $\tQ$ is $\delta$-hyperbolic
(in the sense of Gromov), hence
the definition of $\p_\infty\tilde{Q}$ remains valid, that is the definition using
quasi-geodesics. And this definition coincides with the definition using
geodesics. The most important change to be made here
is a new version of fact 9:\\

\noindent {\bf Lemma 2.5.} {\it The map $\tA:\tN\ra\bo\tQ\sm\bo\tS$, given by $\tA(v)=[c_v]$ is cellular.
Furthermore, we can extend $\tA$ to $\tW \ra\ctQ\sm\bo\tS$ by  defining $\tA(v)=\tilde{E}(\varsigma
(|v|)\frac{v}{|v|})=exp_q (\varsigma  (|v|)\frac{v}{|v|})$, for $|v|<1$, $v\in\tW_q$. 
This extension is cellular and  a homeomorphism on $\tW\sm\tN$.}\\

\noindent {\bf Proof.} Standard Hadamard manifold techniques show that the map $\tA:\tN\ra\bo\tQ\sm\bo\tS$
is continuous.
Let $v\in \tN_p$, $p\in\tS$. We now prove that $C=\tA^{-1}(\tA(v))$ is homeomorphic to a convex
set in $\tS$. First note that for $v, v'\in \tN_p$ we have $\tA(v)\neq \tA(v')$ because 
$c_v$ and $c_{v'}$ are two geodesic rays emanating from the same point. Hence the continuous map $\pi|_{C}: C\ra\pi(C)$
is injective. (Recall $\pi:\tN\ra\tS$ is the bundle projection).\\

Let  $v'\in \tN_{p'}$, $p'\neq p$,  be such that $\tA(v')=\tA(v)$. Let $[p,p']$ be the unique
geodesic segment joining $p$ to $p'$. Since the geodesic rays $c_v$, $c_{v'}$ make a right angle
with $\tS$, at $p$ and $p'$ respectively, we have that $c_v$, $c_{v'}$ and $[p,p']$ bound
a flat geodesic ribbon. Hence $[p,p']\sbs \pi(C)$ and it follows that $\pi (C) $ is convex.
For each $q\in\pi(C)$ there is a unique $v_q\in \tN_q\cap C$ and it is straightforward 
(using the ribbon property) to prove that $q\mapsto v_q$ is continuous. Hence $C$ is homeomorphic
to the convex set $\pi(C)$. Note that $\pi(C)$ is bounded: otherwise $\tQ$ would contain flat
geodesic ribbons isometric to $[0,\ell]\times [0,\infty)$ with $\ell\ra \infty$ which would
contradict the $\delta$-thinness of triangles in the $\delta$-hyperbolic space $\tQ$.
It follows that $C$ is compact. This proves the Lemma.\\

\noindent {\bf Lemma 2.6.} {\it The map $\tA$ is onto and proper.}\\

\noindent {\bf Proof.} Let $p\in \tS$ and $\alpha$ a geodesic ray emanating
from $p$ and not contained in $\tS$. We have to show that there is
a $c_v$ (emanating from some $q\in \tS$ with direction $v$ perpendicular to
$\tS$) such that $\alpha$ and $c_v$ determine the same point in $\tQ$.
Let $c_{v_n}$, at $q_n$, be perpendicular to $\tS$ and passing through
$\alpha(n)$. Using the $\delta$-hyperbolicity of $\tQ$ we get that
$\{ q_n\}$ is bounded so we can assume $q_n\ra q$. Furthermore
we can assume $v_n\ra v$. It is straightforward to verify that
$c_v$ and $\alpha$ determine the same point at infinity. This proves
that $\tA$ is onto.\\

We identify $\p\tQ$ with the unit sphere $\bS$ in $T_p\tQ$, for some $p\in\tS$.
Let $K\sbs \bS-T_p\tS$ be compact. We now prove $\tA^{-1}(K)$ bounded.
Note that if $\alpha$ is a ray emanating from $p$ with direction $\alpha'(0)\in K$, then the angle between $\alpha'(0)$ and $\tS$ is bounded away from zero,
hence there is a $\kappa>0$ such that $d(\alpha(t),\tS)\geq \kappa t$,
for every such $\alpha$. Let $v\in \tA^{-1}(K)$, with $v\in T_q\tS$. Then $[c_v]=[\alpha]$, for some $\alpha$ as above. 
Using the $\delta$-hyperbolicity of $\tQ$ we get
that $d(p,q)$ cannot be arbitrarily large. This proves that
 $\tA^{-1}(K)$ is  bounded. This proves the lemma.\\

\noindent {\bf Lemma 2.7.} {\it The injectivity radius at $p\in Q$ tends to infinity, as
$p$ gets far from $S$.}\\

\noindent {\bf Proof.} Suppose not. Let $\gamma_n$ be non-contractible loops
in $S$ with $d(\gamma_n,S)= n$ and the lengths $\ell(\gamma_n)$ bounded 
(say by $a>0$). Each $\gamma_n$ is homotopic to a closed geodesic
$\beta_n$ in $S$. Lifting to $\tQ$ we obtain $p_n, p'_n\in \tS$ and vectors $v_n,
v'_n$ such that $d(c_{v_n}(n),c_{v'_n}(n))\leq a$ and $b\leq d(p_n,p'_n)\leq a$,
where $b>0$ is the injectivity radius of $S$.
Since $\tS$ has a compact fundamental domain we can assume $p_n\ra p$, $p'_n\ra p'$,
$v_n\ra v$ and $v'_n\ra v'$. It is straightforward to check that
$c_v$ and $c_{v'}$ bound an infinite flat (half) ribbon. We can repeat this
process with the set $\{\gamma_n^k\}$, where $\alpha^k=\alpha*...*\alpha$
(concatenation $k$ times). In this way we get that $\tQ$ contains flat
geodesic ribbons isometric to $[0,b]\times [0,\infty)$ with $b\ra \infty$ which would contradict the $\delta$-thinness of triangles in the $\delta$-hyperbolic space $\tQ$.
This proves the Lemma.\\

Lemma 2.2 of \cite{FO6} remains true. In the proof Lemma 2.6 above is used and there is no need for some
of the (now not valid) facts mentioned in section 2 of \cite{FO6}.  We can now descend to $Q$
and define, as in Lemma 2.4 of \cite{FO6}, the map $A$:\\

\noindent {\bf Lemma 2.8.} {\it
The map $A:N\ra\bo Q$, given by $A(v)=[c_v]$ is cellular.
Furthermore, we can extend $A$ to  $W \ra\bo Q\cup Q$ by  defining 
$A(v) = E( (\varsigma  (|v|)\frac{v}{|v|}))$, for $|v|<1$. This extension is
cellular and a homeomorphism on $W\sm N$.} \\

\noindent {\bf Proof.} Since $\tA$ covers $A$, it is enough to prove that $\gamma C\cap C=\emptyset$,
for $\gamma\in \Gamma$ and $C$ as in the proof of Lemma 2.5.
Suppose there is $v\in \tN_p$ with $\tA(v)=\tA(\gamma_*(v))$, for some nontrivial  $\gamma\in \Gamma$.
Note that $\gamma_*(v)\in\tN_{\gamma(p)}$ and $\gamma(p)\neq p$.
But then  $\tA(v)=\tA(\gamma^n_*(v))$,  $\gamma^n_*(v)\in\tN_{\gamma^n(p)}$ and 
the distance between $\gamma^n(p)$ and $p$ becomes large. Therefore we again obtain large
geodesic flat ribbons in $\tQ$, which cannot happen. This proves the Lemma.\\
\vspace{.6in}

\noindent {\bf \large Changes to Section 3 of \cite{FO6} and Proof of Theorem 2.4.}\\

\noindent {\bf 1.} In the first three paragraphs  replace $\mo(M)$ by $\moo(M)$
and delete the second sentence of the fourth paragraph.\\

\noindent {\bf 2.} The remark at the beginning of section 3 will not be needed because the 
embedding results 1.1-1.4 in section 1 of \cite{FO6} will not be used directly here. 
A new argument is given in item 7 below that does not use (directly) these embedding results
(even though they are used indirectly: the fiber of the map $Emb(\bS^1,Q)\hookrightarrow C^\infty(\bS^1,Q)$ is $(n-5)$-connected, which follows from 1.4 of \cite{FO6}).\\

\noindent {\bf 3.} On item (vi) take $c_2=0$.\\

\noindent {\bf 4.} In the paragraph after item (vi) {\it choose} $\alpha_0$ to be {\it one} $\sigma(u_0)$-geodesic homotopic to $\alpha$.
Also, for $u\in \bS^k$,  {\it define} $\alpha_u=\phi_u(\alpha_0)$. By applying lemma 1.4 of \cite{FO6} (write $\bS^k=\D^k/\p \D^k$) we can assume $\alpha_u=\alpha_0=\alpha$. Hence, for every $u\in\bS^k$, $\alpha$ is a $\sigma(u)$-geodesic.\\  

\noindent {\bf 5.} The function $h$ defined after item (iv') cannot be defined in our case, but facts 1 and 2 still hold.\\

\noindent {\bf 6.} The proof now continues as follows. Let $X$ be the 
cellular disc constructed in section 1E. 
Since $\sigma$ is constant mod $DIFF(M)$ on $\p\D^k$ we also get that $(X,\bS^k)$ is
a cellular $(k+1)$-disc pair (see paragraph before 1.23).
Recall that the elements of $X$ are pairs $(u,\beta)$, $u\in\D^{k+1}$ and $\beta$ is a $\sigma(u)$-geodesic.
We consider $\bS^k\sbs X$ by identifying $u\in\bS^k=\p\D^{k+1}$ with $(u,\alpha)$. Now, the function $h$ in \cite{FO6} is replaced here by
$h:X\ra Emb(\bS^1,Q)\sbs \Omega=C^\infty(\bS^1,Q)$, $h(u,\beta)=\beta$. Note that $h(\bS^k)=\{\alpha\}$.\\

\noindent {\bf 7.} After item (vi) the proof continues as follows.
Since $\pi_k(C^\infty(\bS^1,Q))=0$, $k>1$, and the fiber of $Emb(\bS^1,Q)\hookrightarrow C^\infty(\bS^1,Q)$ is $(n-5)$-connected (this follows from Lemma 1.4 of \cite{FO6}) we have also that $\pi_k(Emb(\bS^1,Q))=0$. This together with Proposition 1.1 and the fact $Emb(\bS^1, Q)$
is locally contractible imply that
$h:(X,\bS^k)\ra (Emb(\bS^1,Q),\alpha)$ is relative null-homotopic (i.e the homotopy
always sends $\bS^k$ to $\{\alpha\}$).
Hence for each $x=(u,\beta)\in X$ there is an isotopy $h_t(x):\bS^1\ra Q$ such that $h_0(x)=\beta$ and $h_1(x)=\alpha$.\\

\noindent {\bf Remark.} By modifying $X$ a little bit we can also assume the above to be true for $k=0$: consider the map
$X\ra Q$, $(u,\beta)\mapsto \beta(1)$. If this map represents $n\in\Z=\check{H}_1 (\bS^1)=\check{H}_1(Q)$ then just replace
$X$ by $\{ (u, e^{-n\pi \eta(u)i}\beta)\,\,|\,\, (u,\beta)\in X  \}$.\\

Write $\sigma(u,\beta)=\sigma (u)$ (we use the same letter). Hence we are now considering $X$ to be the domain of 
$\sigma$. Now claim 1 of section 3 of \cite{FO6} makes sense. The proof is done by extending the isotopies $h_t(x)$
to ambient isotopies as in Lemma 1.4 of \cite{FO6}, and taking the pullback of $\sigma(u)$ by the end of these isotopies.\\

To complete the proof replace  $\D^{k+1}$ by $X$ and $u$ by $x=(u,\beta)$ throughout. Using the same procedure used in \cite{FO6} to construct, for each $u\in\D^{k+1}$, an
$f_u\in P(\bS^1\times\bS^{n-2})$,  we can construct a $f_x\in  CELL(\bS^1\times\bS^{n-2}\times [0,1])$, for $x\in X$. And claim 6 of section 3 of \cite{FO6} still holds with $\iota'\iota\theta$ instead of $\iota\theta$.  By Proposition 1.2
$\iota'\iota\theta$ extends to $\D^{k+1}$. This proves Theorem 2.4.\\

\noindent {\bf Remark.} We have to be careful here since the original
metric is not necessarily negatively curved, hence it may not determine
the identity on $\bS^1\times\bS^{n-2}\times [0,\infty]$. To overcome
this let $B_t:W\ra\bo Q\cup Q$ be a 1-parameter family of cellular maps
such that $B_0$ is a homeomorphism and $B_1$ is the map $A$ of lemma 2.8
corresponding to the base point $u_0\in \bS^k$. Then identify
$W=\bS^1\times\bS^{n-2}\times [0,\infty]$ with $\bo Q\cup Q$
by $B_0$.

\vspace{.8in}

\noindent {\bf \Large  Appendix A. The Swinging Neck.}\\

For a function $h:\R\ra (0,\infty)$ denote by $M_h$ the surface of revolution obtained by rotating
the graph  $\{ (x,h(x), 0 )\, :\,\, x\in \R\}$ of $h$ around the $x$-axis.
We consider $M_h$ with the Riemannian metric induced by $\R^3$.\\

Let $f:\R\ra [1,\infty)$ be a smooth function such that: {\bf (1)} $f\equiv 1 $ on $[-1,1]$, \,\,\,
{\bf (2)} $f''(x)> 0$, $| x|> 1$,\,\,\, {\bf (3)} $f''(x)\geq \delta >0$, for $|x|\geq 2$.
Then $M_f$ is nonpositively curved and contains the flat cylinder $[-1,1]\times\bS^1$.\\

Let $\alpha:\R\times[-2,2]\ra [0,1]$ be a smooth function such that (we write $\alpha_t$ for the function $x\mapsto \alpha(x,t)$):
{\bf (1)} $\alpha\equiv 0 $ for  $|x|\geq 4$ and all $t$, \,\,\, {\bf (2)}  $\alpha_t''(x)>0$, for $|x|\leq 3$ and all $t$, \,\,\,
{\bf (3)} $\alpha_t$ has a unique minimum value (equal to 0)  on $[-3,3]$ at $t$, for all $t$.\\

Define $F:\R\times [0,1]\ra [1,\infty)$ by $F(x,0)=f(x)$, and for $t\in (0,1]$ by

$$ F(x,t)=f(x)+e^{-1/t}\alpha (\, x, \sin\, (1/t)\, )$$

\noindent and write $f_t(x)=F(x,t)$. Thus $f_0=f$. Then $F$ is smooth and for small enough $t>0$ we have:
{\bf (1)}  $f_t''(x)\geq 0$, $\forall x\in \R$\,\,\, 
{\bf (2)} $f_t$ has a unique minimum value at $\sin \,(1/t)$.
{\bf (3)} $f_t\equiv f$ outside $[-4,4]$.\\

Write $M_t=M_{f_t}$ and $M=M_f$. Then $M_t$  is negatively curved and concides with $M$ outside a compact set.
Note that all  $\{ x\}\times\bS^1$, $x\in [-1,1]$, are non-trivial closed geodesics of minimal length in  $M$.
But $M_t$ has a unique non-trivial closed geodesic $\{\sin\, (1/t) \}\times\bS^1$ of minimal length that oscillates 
between $\{ -1\}\times\bS^1$ and $\{ 1\}\times\bS^1$ faster and faster, as $t$ approaches 0.\\

Note that, with some care,  we can fit these ``necks''- the relevant parts of $M_t$ and $M$- on a closed negatively curved surface.

F.T. Farrell

SUNY, Binghamton, N.Y., 13902, U.S.A.\\

P. Ontaneda

SUNY, Binghamton, N.Y., 13902, U.S.A.

\end{document}